\documentclass[a4paper, 12pt]{article}
\usepackage[utf8]{inputenc} 
\usepackage{makeidx}
\usepackage{enumerate, theorem}
\usepackage{amsmath, amsfonts, amssymb}
\usepackage[height=22.5cm, width=15cm]{geometry}
\usepackage[all]{xy}
\usepackage{mathrsfs, yfonts}
\usepackage{graphicx}
\newcommand{\parag}[1]{\paragraph{\sc{#1.}} }
\newtheorem{thm}{Theorem}[subsection]
\newtheorem{defn}[thm]{Definition}
\newtheorem{cor}[thm]{Corollary}
\newtheorem{prop}[thm]{Proposition}
\newtheorem{lemma}[thm]{Lemma}

\begin{document}
\title{ On Principal Value and Standard Extension of Distributions}

 \author{Daniel Barlet\footnote{Institut Elie Cartan, G\'eom\`{e}trie,\newline
Universit\'e de Lorraine, CNRS UMR 7502   and  Institut Universitaire de France.}.}

\maketitle

\hfill \qquad  \qquad \qquad {\it To the memory of Jan Erik Bj\"{o}rk who }

\hfill \qquad \qquad \qquad   {\it explains to me the theory of D-modules }

\parag{Abstract}  For a holomorphic function $f$ on a complex manifold $\mathscr{M}$ we explain in this article that the distribution associated to $\vert f\vert^{2\alpha}(Log\vert f\vert^2)^q f^{-N}$ by taking the corresponding limit on the sets $\{ \vert f\vert \geq \varepsilon \}$ when $\varepsilon $ goes to $0$, coincides for $\Re(\alpha) $ non negative and $q, N \in \mathbb{N}$, with the value at $\lambda = \alpha$ of the meromorphic extension of the distribution $\vert f\vert^{2\lambda}(Log\vert f\vert^2)^qf^{-N}$. This implies that any distribution in  the $\mathcal{D}_{\mathscr{M}}$-module generated by such a distribution has the Standard Extension Property. This implies a non torsion result for the $\mathcal{D}_{\mathscr{M}}$-module generated by such a distribution. As an application of this result we determine  generators for the conjugate modules of the regular holonomic $\mathcal{D}$-modules associated to $z(\sigma)^\lambda$, the power $\lambda $, where $\lambda$ is any complex number,  of the (multivalued) root of the universal equation of degree $k$, $z^k + \sum_{j=1}^k (-1)^h\sigma_hz^{k-h} = 0$ whose structure is studied in \cite{B4}.
 
 \parag{AMS Classification}32 A 27- 32 D 15- 32 S 40.

\parag{Key words} Principal Value- Meromorphic Extension- Regular Holonomic Distribution- Kashiwara Conjugation Functor- Standard Extension Property for Distributions.

\tableofcontents

\section{Introduction}

Let $f : \mathscr{M} \to \mathbb{C}$ be a holomorphic function on a domain $\mathscr{M} $ in $\mathbb{C}^{n+1}$  and assume that $f$ has no critical value different from $0$. 
Let $\alpha$ be a  complex number with a non negative real part  and $N$ a positive integer. There are two natural ways  to define a distribution on $\mathscr{M} $ whose restriction to $ \mathscr{M}  \setminus \{f = 0\}$ is equal to $\vert f\vert^{2\alpha}f^{-N}$. The first one is given by the principal value method (see for instance \cite{H-L}):\\

Let $\xi$ be a $\mathscr{C}^\infty_c(\mathscr{M} )$ differential form on $\mathscr{M} $ of type $(n+1, n+1)$ and define for $\varepsilon > 0$ the distribution $T_{\alpha,N}^\varepsilon $ by
$$ \langle T_{\alpha,N}^\varepsilon, \xi \rangle := \int_{\vert f\vert \geq \varepsilon} \vert f\vert^{2\alpha}f^{-N}\xi.$$
 We  give in section 2 a rather short proof of the fact that the limit  $T_{\alpha,N}$ when $\varepsilon$ goes to $0$ of $ \langle T_{\alpha,N}^\varepsilon, \xi \rangle $ exists and defines a distribution $T_{\alpha,N}$ on $\mathscr{M} $ whose restriction to $ \mathscr{M}  \setminus \{f = 0\}$ is equal to $\vert f\vert^{2\alpha} f^{-N}$. \\

The second method is to show that for $\Re(\lambda)$ large enough and for  any test differential  form  $\xi \in \mathscr{C}^\infty_c(\mathscr{M} )^{(n+1,n+1)}$  the function
$$ \lambda \mapsto \int_{\mathscr{M} } \vert f\vert^{2\lambda}f^{-N}\xi$$
is holomorphic and defines a holomorphic family of distributions on $\mathscr{M} $. Moreover this holomorphic family of  distributions admits a meromorphic extension to the complex $\lambda$-plane with no pole when $\Re(\lambda) \geq 0$. Then,  
$\alpha $ is not a pole and  we define the distribution $S_{\alpha,N}$ on $\mathscr{M} $ as the value at $\lambda = \alpha$ of this meromorphic extension.\\

We give a proof of the existence of $T_{\alpha,N}$ and $S_{\alpha,N}$ and we prove  the equality  in $Db_{\mathscr{M}}$, the sheaf of distributions on $\mathscr{M}$, $T_{\alpha,N} = S_{\alpha,N}$.\\
 Our main tools for these proofs are 
\begin{enumerate}
\item The Asymptotic Expansion Theorem of fiber-integrals given in \cite{B2}  (which uses Hironaka's Desingularization Theorem); see also \cite{BMa}.
\item The existence of a local Bernstein identity for $f$, due to J.E. Bjork in the analytic case (see \cite{Be} and  \cite{Bj1}).
\end{enumerate}

As an application we first deduce of this result the absence of torsion for some $\mathcal{D}$-modules generated by the distributions constructed in the first part. In fact we prove more : any distribution in such a $\mathcal{D}$-module has the Standard Extension Property (compare with the result of \cite{B-K}). We recall the  definition of the Standard Extension Property in the begining of the section 2.\\
Then we apply this Theorem to the determination of generators of  the conjugate $\mathcal{D}$-modules of the $\mathcal{D}$-modules $\mathcal{N}_\lambda$  associated to $z(\sigma)^\lambda$, the power $\lambda $, where $\lambda$ is any complex number, of the (multivalued) root of the universal equation of degree $k$, $z^k + \sum_{j=1}^k (-1)^h\sigma_hz^{k-h} = 0$ whose structure is studied in \cite{B4}.\\

\section{Existence of the Principal Value}

\subsection{The Standard Extension Property}

First, we recall the definition of this property (see  \cite{B3} in  Appendix and \cite{B-K} in Paragraph 3.1 or \cite{Bj2}) 

\begin{defn} \label{SEP}
Let $\mathscr{M} $ be a complex manifold of pure dimension $n+1$ and let  $T$ be a distribution on $\mathscr{M} $. We shall say that $T$ has the {\bf Standard Extension Property} if the following condition are fulfilled:
\begin{enumerate}
\item Outside a hypersurface $H$ in $\mathscr{M} $ the distribution  $T$ is a $\mathscr{C}^\infty$ function.
\item For each point $x$  in $H$ there exists an open neighborhood $U$ of $x$, a local holomorphic  equation $\{f = 0\} = H \cap U$ of $H$ in $U$ such that for any test differential form $\xi$  in $\mathscr{C}_c^\infty(U)^{(n+1,n+1)}$ we have
$$ \langle T, \xi \rangle = \lim_{\varepsilon \to 0} \int_{\vert f\vert \geq \varepsilon}  T\xi .$$
\end{enumerate}
\end{defn}

It is easy to see that if the condition 2 of the previous definition is satisfied for some choice of local equation of $H$, it is satisfied for any other choice.\\ 
Note also that this property is clearly stable by multiplication by a $\mathscr{C}^\infty$ function but is not stable, in general, by the action of $\mathcal{D}_{\mathscr{M}}$. For instance the locally integrable function  $1/\bar z$ defines a distribution on $\mathbb{C}$  which has the standard extension property, but $\partial_z(1/\bar z) = i\pi\delta_0$ is a non zero torsion element in $Db_{\mathbb{C}}$.\\

The following result is proved in \cite{B-K}

\begin{thm}\label{B-K}
Let $\mathcal{M}$ be a  regular holonomic $\mathcal{D}_{\mathscr{M}}$-module and let  $H$ be a hypersurface in $\mathscr{M} $ such that $\mathcal{M}_{\mathscr{M}  \setminus H}$ is $\mathcal{O}_{\mathscr{M}  \setminus H}$-coherent. Let $i : \mathcal{M} \to Db_{\mathscr{M}}$
be a $\mathcal{D}_{\mathscr{M}}$-linear morphism such that $i $ belongs to $L^2(c(\mathcal{M}))$, the $L^2$ lattice (see \cite{B-K})  of the conjugate module of $\mathcal{M}$ (see \cite{Ka1}). Then any distribution in $i(\mathcal{M})$ has the Standard Extension Property (relatively to $H$).
\end{thm}

The main difficulty to use this theorem is the identification of the $L^2$-lattice. When the $\mathcal{D}_{\mathscr{M}}$-module has no $\mathcal{O}_{\mathscr{M}}$-torsion, the $L^2$-lattice is given by a  local square integrability condition along $H$ which is often rather easy to verify. But we must know "a priori" that the $\mathcal{D}_{\mathscr{M}}$-module we are considering has no $\mathcal{O}_{\mathscr{M}}$-torsion to use this simple characterization of the $L^2$-lattice. We shall  see in section 5 that the main point is to prove that there  is no torsion in some of the regular holonomic $\mathcal{D}_{\mathscr{M}}$-modules involved in our computations. 

\bigskip

\subsection{Principal Value}

Let us begin by giving the precise Asymptotic Expansion Theorem for fiber-integrals proved in \cite{B1} which we shall need below. We keep the situation and  the notations introduced in the beginning of section 1.

\begin{thm}\label{AS.EX.1}
Let $\varphi \in \mathscr{C}^\infty_c(\mathscr{M} )$ be a differential form of type $(n, n)$ on $\mathscr{M} $. Define the function $\Theta(s) := \int_{f=s} \varphi $ for each $s \in \mathbb{C}$. This function, which is $\mathscr{C}^\infty$ on $\mathbb{C} \setminus \{0\}$, admits when $s \to 0$ a asymptotic expansion of the form
\begin{equation}
\Theta(s) \simeq \sum_{m,m'}^{r, j}  a_{m,m'}^{r,j}\vert s\vert^{2r}(Log\, \vert s\vert)^js^m\bar s^{m'}
\end{equation}
where $m, m'$ are non negative integers, $r$ describes a finite set $R \subset [0, 1[ \cap \mathbb{Q}$ and where $j$ is an integer in $[0, n+1]$. The finite set $R$ is independent of the choice of  $\varphi$.  Moreover, this asymptotic expansion is term-wise differentiable at any order.$\hfill \blacksquare$\\
\end{thm}

\parag{Remarks}\begin{enumerate}
\item The continuity of the function $\Theta$ at the point $s = 0$ implies that for $r = 0$ and $j \geq 1$ we have $m+m' \geq 1$. In fact it is proved in Proposition 6 of {\it loc. cit.} that for $j \geq 1$ we have $a_{m, 0}^{0, j} = 0$ (resp. $a_{0,m'}^{0, j} = 0$) which implies that for $r= 0$ and $j \geq 1$ we have $m \geq 1$ and $m' \geq 1$ when $a_{m,m'}^{0,j} \not= 0$. This shows that no new type of term appears in such an expansion when we apply $s\partial_s$ or $\bar s\partial_{\bar s}$. And, of course, neither a constant term nor  a term like  $s^m(Log\vert s\vert)^j$ or $\bar s^{m'}(Log\vert s\vert)^j$  with $j \geq 1$ may appear in these expansions. So the functions $s\partial_s\Theta, \bar s\partial_{\bar s}\Theta$ and $s\bar s\partial_s\partial_{\bar s}\Theta$  are bounded by $O(\vert s\vert^\gamma)$ for some $\gamma > 0$ when $s$ goes to $0$.
\item It is also proved, in {\it loc. cit.} that the linear map $\varphi \mapsto a_{m,m'}^{r,j}(\varphi)$ is a $(1,1)$-current on $\mathscr{M} $. This current is supported by the singular set of the  hypersurface $\{f = 0\}$ when $(r, j) \not= (0, 0)$, that is to say for terms which do not correspond to an usual term of a Taylor expansion at the origin of a $\mathscr{C}^\infty$ function.\\
\end{enumerate}

It will be important in the sequel to use the following corollary of this theorem.

\begin{cor}\label{AS.EX.2}
In the situation of the previous theorem, for any differential $\mathscr{C}^\infty$ form $\psi$ and $\xi$ respectively  of type $(n,n+1)$ and $(n+1,n+1)$ with $f$-proper supports in $\mathscr{M}$ the fiber-integrals
$$ \eta(s) := \int_{f=s} \bar f\frac{\psi}{d\bar f} = \bar s \int_{f=s}  \frac{\psi}{d\bar f} \quad {\rm and} \quad \zeta(s) := \int_{f=s} f\bar f\frac{\xi}{df\wedge d\bar f} = s\bar s\int_{f=s} \frac{\xi}{df\wedge d\bar f} $$
admit asymptotic expansions of the same type as above when $s$ goes to $0$ (but see the previous remark 1). 
\end{cor}

\parag{Remark} For any compact set $K$ in $\{f = 0\}$ there exists an integer $N$ such that we have the inclusion of sheaves in an open neighborhood of $K$ in $\mathscr{M} $
$$ f^N\Omega^{n+1}_\mathscr{M}  \subset df \wedge \Omega^n_\mathscr{M}  \quad {\rm and} \quad f^N\bar f^N \Omega^{n+1}_\mathscr{M} \wedge \overline{\Omega}^{n+1}_\mathscr{M}  \subset \Omega^n_\mathscr{M} \wedge \overline{\Omega}^n_\mathscr{M} \wedge df\wedge \bar df .$$
So the only new information in the previous corollary is that we may keep $m$ and $m'$ non negative in the expansions.

\parag{Proof} Remark first that the previous theorem is in fact local around the hypersurface $\{f = 0\}$ in $\mathscr{M}$ and so the asymptotic expansion is valid for any $\mathscr{C}^\infty$  differential form of type $(n,n)$ with f-proper support in $\mathscr{M} $.\\
Then using the fact that for such a form $\varphi$ we have the same type of asymptotic expansion for
 $$s\partial_s\Theta(s) =  \int_{f=s} \frac{fd'\varphi}{df}, \quad  \bar s\partial_{\bar s}\Theta(s) =   \int_{f=s} \frac{\bar fd''\varphi}{d\bar f} \quad {\rm and\ for} \quad  s\bar s\partial_s\partial_{\bar s} \Theta(s) = \int_{f=s} \frac{f\bar fd'd''\varphi}{df\wedge d\bar f} $$
 we see that it is enough to show that we have the local surjectivity of $d', d'', d'd''$ with f-proper supports   for the type $(n+1, n), (n, n+1)$ and $(n+1, n+1)$ respectively.\\
 Using a partition of unity, this is consequence of the following lemma and the local parametrization theorem for the hypersurface $\{f = 0 \}$\footnote{In order that we may choose, near each point in $\{ f = 0 \}$, $U$ and $D$ such that   $\bar U \times \partial D$ does not meet $\{f= 0\}$ in our next lemma.}.
 
 \begin{lemma}\label{Dolbeault}
 Let $U$ be an open polydisc in $\mathbb{C}^n$ and $D$ a disc in $\mathbb{C}$. Let $\psi$ and $\xi$ be $\mathscr{C}^\infty$ forms respectively of type $(n,n+1)$ and $(n+1, n+1)$ with support in $K \times D$ where $K \subset U$ is a compact set. Fix a relatively compact open disc $D'$ in $D$. Then there exists $\varphi_1$ and $\varphi_2$ which are $\mathscr{C}^\infty$ forms of type $(n,n)$ and with support in $K \times D$ such that
 $$ d''\varphi_1 = \psi \quad {\rm and} \quad d'd''\varphi_2 = \xi $$
 on $U\times D'$.$\hfill \blacksquare$
 \end{lemma}
 
 \parag{Proof} This an easy consequence of the fact that $\partial_z$ and $\partial_z\partial_{\bar z}$ are surjective on $\mathscr{C}^\infty(D)$ and that we may solve the corresponding equations with $\mathscr{C}^\infty$ dependence of a parameter on a relatively compact open disc $D' \subset\subset D$ using a fundamental solution of the corresponding operator ( see for instance  \cite{B-M} chapter IV Proposition 5.2.4 for details). $\hfill \blacksquare$
 
\parag{Remark}  Note that $d'\bar \varphi_1 = \bar\psi$ gives also the $d'-$case.

\begin{thm}\label{PV.1}
In the situation $f : \mathscr{M}  \to \mathbb{C}$ introduced in section 1, let $\alpha$ be a complex number such that its real part $\Re(\alpha)$ is non negative. Then, for any positive integer $N$ and  for any $\mathscr{C}^\infty$ differential form  $\xi$ of type $(n+1,n+1)$ with compact support in $\mathscr{M} $ the limit when $\varepsilon > 0$ goes to $0$ of 
$$  \langle T_{\alpha,N}^\varepsilon, \xi \rangle := \int_{\vert f\vert \geq \varepsilon} \vert f\vert^{2\alpha}f^{-N}\xi$$
exists and defines a distribution (that is to say a $(0,0)$ current)  $T_{\alpha,N}$ on $\mathscr{M} $. 
\end{thm}

The main argument to prove this result uses  the following consequence of the previous corollary of the Asymptotic Expansion Theorem (see also \cite{H-L} for a proof of the proposition below using a direct computation in a desingularization of the hypersurface $\{f = 0 \}$).

\begin{prop}\label{P.V.2}
Let $\psi$ be a $\mathscr{C}^\infty$ differential form of type $(n, n+1)$ with f-proper support in $\mathscr{M} $. Then we have for any  positive integer $N$
$$ \lim_{\varepsilon \to 0} \   \int_{\vert f\vert = \varepsilon} f^{-N}\psi= 0 .$$
\end{prop}

\parag{Proof} Note first that the result is local and is obvious near a point where $f$ does not vanish. Around a point where $f$ vanishes,  Milnor's Fibration Theorem  allows to use Fubini's Theorem to compute this integral for $\varepsilon$ small enough  as follows:
$$ \int_{\vert f\vert = \varepsilon} f^{-N}\psi = \int_{0}^{2\pi} d\theta \int_{f = s} f^{-N}\frac{\psi}{d\theta} $$
where $\psi /d\theta$ on $\{f = s = \varepsilon  e^{i\theta} \}$ is equal to \  $\bar f\psi/id\bar f$ \ because we have $d\theta = id\bar f/\bar f$ and taking in account the type of $\psi$ and the fact that $\{f = s\}$ is a complex $n$-dimensional sub-manifold in $\mathscr{M}$ for $s \not= 0$. This gives with $s := \varepsilon e^{i\theta}$
$$ \int_{\vert f\vert = \varepsilon} f^{-N}\psi = \int_{0}^{2\pi} d\theta \int_{f = s} f^{-N}\frac{\bar f\psi}{id\bar f}. $$
When the integral $\int_{f = s} \bar f\psi/d\bar f$ is $\mathcal{O}(\varepsilon^{N+1})$ it is clear that the limit when $\varepsilon$ goes to $0$ exists and vanishes. So it is enough to show the result when we replace the fiber-integral 
$$-i\ s^{-N}\bar s\partial_{\bar s}\Theta = \int_{f=s} f^{-N} \bar f\frac{\psi}{id\bar f}$$ 
 by its asymptotic expansion at a high enough order, and then, by linearity, to prove the assertion when we replace the fiber-integral by $ \vert s\vert^{2r}(Log \vert s\vert)^js^{m-N}\bar s^{m'} $. In this case the exponent of $e^{i\theta}$ is given, for $s = \varepsilon e^{i\theta}$, by :
$$ -N + m - m'$$
and the corresponding exponent of $\varepsilon$ is given by 
$$ 2r +m +m' - N.$$
In order to find a non zero integral between $0$ and $2\pi$ we need that $m = m' + N$ and then the exponent of $\varepsilon$ is then equal to $2r + 2m'$. So the only terms  where the limit is not clearly equal to $0$ appear when   $r = m' = 0$. But in this case the limit  is again zero thanks to Remark 1 following Theorem  \ref{AS.EX.1}.$\hfill \blacksquare$

\parag{Proof  of Theorem \ref{PV.1}} The assertion is local near each point of $\{f = 0\}$ so, using Lemma \ref{Dolbeault} (in fact the remark following it) we may assume that $\xi = d'\psi$ where $\psi$ is a $\mathscr{C}^\infty$ differential form of type $(n, n+1)$ with f-proper support. Then Stokes' Formula gives
$$(\alpha - N) \int_{\vert f\vert \geq \varepsilon} \vert f\vert^{2\alpha}f^{-N}\frac{df\wedge\psi}{f} +  \int_{\vert f\vert \geq \varepsilon} \vert f\vert^{2\alpha}f^{-N}\xi = \int_{\vert f\vert = \varepsilon} \vert f\vert^{2\alpha}f^{-N}\psi.$$
We know that the limit of the right hand-side is $0$ when $\varepsilon$ goes to $0$ thanks to Proposition \ref{P.V.2}.
Now write the first integral  in the left hand-side as
$$ \int_{\vert s\vert \geq \varepsilon} \vert s\vert ^{2\alpha}s^{-N} \frac{ds\wedge d\bar s}{s\bar s} \int_{f= s} \frac{\bar f\psi}{d\bar f} $$
and using polar coordinates $s = \rho e^{i\theta}$ so $ds \wedge d\bar s = -2i \rho d\rho d\theta$ and the asymptotic expansion of the fiber-integral of  $\bar f \psi / d\bar f $ we may replace this fiber-integral by its asymptotic expansion at the origin at a sufficient large order. Then, using  linearity, it is enough  to consider the integrability at $0$ of each non zero term which only depends on the real part of $\alpha$. So we have only to consider the terms for which
$$ -N + m-m' = 0 .$$
The corresponding real part of the power of $\rho$ is given by 
$$-N - 1 + 2r + m + m' + 2\Re(\alpha) = -1 + 2r + 2m' + 2\Re(\alpha).$$
 This  is  at least equal to $-1$  for $\Re(\alpha ) \geq 0$ because $m' + r \geq 0$. Then   either we have  $r + m' + \Re(\alpha) > 0$ and this implies the integrability at $0$, or $r = m' = \Re(\alpha) =  0$. And in this case, thanks to  Remark 1 following Theorem \ref{AS.EX.1}, the term to integrate  is $O(\vert s\vert^{\gamma-1})$ for some $\gamma > 0$, giving again the local integrability at $0$.  So the limit of $\langle T_{\alpha,N}^\varepsilon , \xi \rangle$ exists for each  test form $\xi$ on $\mathscr{M} $ when $\varepsilon$ goes to $0$.\\
 The fact that the so obtained linear form on test differential  forms is a distribution is an easy exercise left to the reader.$\hfill \blacksquare$\\
 
 Note that in the previous proof, we may conclude with the weaker hypothesis asking that  $\Re(\alpha) + r > 0$ for any $r$ in  $R \cup \{1\} \setminus \{0\}$. Remark that with this hypothesis $\Re(\alpha)$ is not a pole for the meromorphic extension of $F_{ N,\xi}(\lambda), \forall N \in \mathbb{N}$, defined below.

\parag{Remark $R1$}
If we replace $\vert f\vert^{2\alpha}$ by $\vert f\vert^{2\alpha}(Log\vert f\vert^2)^\beta$ where $\beta$ is any positive real number, the same result holds true with the same proof.

\section{Bernstein identity and meromorphic extension}

We first recall the fundamental theorem of Bernstein \cite{Be}, generalized to the local analytic case by Bjork \cite{Bj1}.

\begin{thm}\label{Bernstein 0}
Let $f : (\mathbb{C}^{n+1}, 0) \to (\mathbb{C}, 0)$ be a non zero holomorphic germ. Then there exists a monic polynomial \ $b \in \mathbb{C}[\lambda]$ and, near $0$ in $\mathbb{C}^{n+1}$,  a holomorphic partial differential operator $P$ with polynomial coefficients in $\lambda$ such that the identity
\begin{equation}
 P(z, \partial_z, \lambda)f^{\lambda+1} = b(\lambda)f^\lambda
 \end{equation}
holds locally outside the hypersurface $\{f = 0 \}$ on some open neighborhood of the origin.
\end{thm}

The minimal monic polynomial $b \in  \mathbb{C}[\lambda]$ such that and identity of this kind holds near $0$  is called the {\bf Bernstein polynomial of $f$ at the origin}.\\

Recall also the fundamental result of Kashiwara \cite{K}.

\begin{thm}\label{Bernstein 1}
For any non zero holomorphic  germ $f$  at the origin, the roots of the Bernstein polynomial of $f$ at the origin are rational and negative.$\hfill \blacksquare$
\end{thm} 

\begin{cor}\label{Bernstein 2}
Let $f : (\mathbb{C}^{n+1}, 0) \to (\mathbb{C}, 0)$ be a non zero holomorphic germ and let $\mathscr{M} $ be an open neighborhood of the origin on which the Bernstein identity of $f$ is valid\footnote{This means precisely that the Bernstein identity for $f$ is valid in the universal cover of \\ $\mathscr{M}  \setminus \{f = 0 \}$, on which $f^\lambda $ is  defined as $\exp(\lambda Log\, f)$ for any given determination of $Log\, f$.}. For any $N \in \mathbb{N}$ and any test differential form  $\xi \in \mathscr{C}^{\infty}_c(\mathscr{M} )^{(n+1,n+1)}$ the function defined for $\lambda \in \mathbb{C}$ such that $\Re(\lambda) \gg N$  by
$$ F_{N, \xi}(\lambda) := \int_\mathscr{M}  \vert f\vert^{2\lambda}f^{-N}\xi$$
is holomorphic and admits a meromorphic extension to the complex $\lambda$-plane with poles of order at most $n+1$ at points in $\cup_r \   \{r - \mathbb{N}\}$ where $r$ is a root of the Bernstein polynomial of $f$ at the origin.\\
Moreover, for any $\alpha \in \mathbb{C}$, the linear forms on $\mathscr{C}^{\infty}_c(\mathscr{M} )^{(n+1,n+1)}$ given by the coefficient $P_k(\lambda = \alpha, F_{N,\xi}(\lambda))$ of $(\lambda - \alpha)^{-k}, k \in \mathbb{Z},  k \leq n+1$, in the Laurent expansion at $\lambda = \alpha$ of the meromorphic extension of $F_{N, \xi}$ is a distribution on $\mathscr{M} $.\\
 For $k \in [1, n+1]$ this distribution has support in $\{f = 0 \}$ and, for $k \in [2, n+1]$ or for $k \in [1, n+1]$ and $\alpha \not\in \mathbb{Z}$, the support of this distribution is contained in the singular set of $\{f = 0 \}$.
\end{cor}

\parag{proof} The equation $(2)$ implies, as $b$ has rational coefficients,  that we may find for any positive integer $M$ a anti-holomorphic differential operator $P_M$ depending polynomially of $\lambda$ such that  we have
\begin{equation}
P_M(\bar f^{\lambda+M}) = b(\lambda)\dots b(\lambda+M-1)\bar f ^{\lambda}.
\end{equation}
Note that the roots of $B_M(\lambda) := b(\lambda)\dots b(\lambda+M-1)$ are negative rational numbers. So for $M \gg N$ we obtain the equality of continuous functions for $\Re(\lambda) \gg N$
$$ \vert f\vert^{2\lambda}f^{-N} = \frac{1}{B_M(\lambda)}P_M(\vert f\vert^{2\lambda}f^{-N}\bar f^M) $$
 If $P_M^*$ is the adjoint of $P_M$ this implies for $\Re(\lambda) \gg N$ and $M$ large enough, that  for any test differential  form  $\xi \in \mathscr{C}^{\infty}_c(\mathscr{M} )^{(n+1,n+1)}$ we obtain
\begin{equation}
 \int_\mathscr{M}  \vert f\vert^{2\lambda}f^{-N}\xi = \frac{1}{B_M(\lambda)} \int_\mathscr{M}  \vert f\vert^{2\lambda}f^{-N}\bar f^MP_M^*(\xi) .
 \end{equation}
This gives the meromorphic extension to the complex $\lambda$-plane of the holomorphic distribution $\xi \mapsto  \int_M \vert f\vert^{2\lambda}f^{-N}\xi $ defined for $\Re(\lambda) \gg N$, because $P_M^*(\xi)$  is $\mathscr{C}^\infty_c$ in $\mathscr{M} $ and depends polynomially on $\lambda$ and because the right hand-side of the above formula is holomorphic on the open set $\Re(\lambda) > -m$ for any given positive integer $m$ as soon as $M$ is large enough compare to $N+m$.\\
Moreover this meromorphic extension has no pole at points  which are not inside the union of the sets $r- \mathbb{N}$ where $r$ is a root of $b$.\\
It is easy to see that near points where $f$ does not vanish, this meromorphic extension has no pole and that near points where $f= 0$ but where $df$ does not vanish the poles of this meromorphic extension are at most simple poles at negative integers. This complete the proof. $\hfill  \blacksquare$\\

\parag{Remark $R2$} Let $q$ be a positive integer. The $q$-th derivative in $\lambda$ of the holomorphic function $F_{N, \xi}(\lambda) $ is given, for $\Re(\lambda) $ large enough by the absolutely converging integral
$$ \int_\mathscr{M}  \vert f\vert^{2\lambda}(Log \vert f\vert^2)^qf^{-N}\xi$$
and the meromorphic extension of these functions allows to define, for each integer $q$, a  meromorphic distributions on $\mathscr{M} $ which has no pole for $\Re(\lambda) \geq 0$.\\
Now  analog arguments as above give easily a generalization of the previous theorem to these cases.

\section{Equality of the Principal Value with the value of the Meromorphic Extension }

\subsection{The equality theorem}

We keep the notations of the introduction.

\begin{defn} In the situation above we define  the distribution $S_{\alpha,N}$ on $\mathscr{M} $ by the formula
$$ \langle S_{\alpha,N} , \xi \rangle := P_0(\lambda = \alpha,  \int_\mathscr{M}  \vert f\vert^{2\lambda}f^{-N}\xi ).$$
\end{defn}

The aim of this  paragraph is to prove the following result:

\begin{thm}\label{coincide 1}
Assume that $\Re(\alpha) \geq 0$. Then for any positive integer $N$ we have for any test differential  form $\xi$:
$$ \langle S_{\alpha,N}, \xi \rangle = \langle T_{\alpha,N}, \xi \rangle =  \lim_{\varepsilon \to 0} \int_{\vert f\vert \geq \varepsilon} \vert f\vert^{2\alpha}f^{-N}\xi .$$
\end{thm}

\parag{Proof} We want to show that the analog of the equality $(4)$ for $\lambda = \alpha$ holds if we perform the integration only on the subset $\{\vert f\vert \geq \varepsilon\}$ with an error which goes to zero when 
$\varepsilon $ goes to $0$. This would be enough to complete the proof. But, of course, the error comes from the boundary terms which are integrals on $\{\vert f \vert = \varepsilon \}$ appearing in the various Stokes Formulas necessary to pass from $P_M$ to its adjoint $P_M^*$. It is easy to see that such "error" terms have the following shape: a polynomial in $\lambda$ with coefficient like
$$ \int_{\vert f\vert = \varepsilon} \vert f\vert^{2\alpha} f^{-N}\bar f^{M'}\psi $$
where $\psi$ is in  $\mathscr{C}^{\infty}_c(\mathscr{M} )^{(n+1,n)}$ and $M'$ is an integer in $[0,M]$. Now using the same arguments than in the proof of Theorem \ref{PV.1} we see that the only non zero term in such an integral comes from the coefficient of $s^m\bar s^{m'}$ in the asymptotic expansion at $s = 0$ of the function $s \mapsto \int_{f=s} \bar f\psi/d\bar f$ such that 
$$ m - m' - N - M' = 0 .$$
And this non zero term comes with a power of $\varepsilon$ which is at least equal to  
$$2\Re(\alpha) + m + m' - N + M' + 2r  = 2\Re(\alpha) + 2m' +2M' + 2r  \geq 0 $$ 
and may be some $(Log\, \varepsilon)^q$ factor. So such term goes to $0$ when $\varepsilon$ goes to $0$ when $m'+r+M' > 0$  for $\Re(\alpha) \geq 0$ but also  in the case where $\Re(\alpha) + m' + M' + r = 0$, thanks again to Remark 1 following Theorem \ref{AS.EX.1}. This  concludes the proof.$\hfill \blacksquare$\\

Again in the previous proof, we may conclude with the weaker hypothesis asking that  $\Re(\alpha) + r > 0$ for any $r$ in  $R \cup \{1\} \setminus \{0\}$.

\parag{Remark $R3$} Using Remarks $R1$ and $R2$ we obtain again with the same proof, that the previous theorem is still valid if we replace $\vert f\vert^{2\alpha}$ by $\vert f\vert^{2\alpha}(Log \vert f\vert^2)^q$ for any positive integer $q$.

\subsection{Non torsion of the corresponding $\mathcal{D}$-modules}

We shall deduce from Theorem \ref{coincide 1} an important  corollary. To formulate this result we need the following definition, where we keep the situation described in the introduction.

\begin{defn}\label{formal action}Let $\alpha$ be a complex number with a non negative real part. Let $V$ be an holomorphic vector field on $\mathscr{M} $. We define the {\bf formal action} of $V$ on $\vert f\vert^{2\alpha} f^{-N}$ by the formula
$$ \langle V, \vert f\vert^{2\alpha} f^{-N} \rangle := (\alpha - N)V(f)\vert f\vert^{2\alpha}f^{-N-1} .$$
Then this defines a "formal action" of $\mathcal{D}_{\mathscr{M}}$ on the $\mathcal{O}_{\mathscr{M}}$-module $\mathcal{O}_{\mathscr{M}}\vert f\vert^{2\alpha}[ f^{-1}]$.
\end{defn}

Remark that, thanks to Theorem \ref{coincide 1}, each element in  $\mathcal{O}_{\mathscr{M}}\vert f\vert^{2\alpha}[ f^{-1}]$ defines an unique distribution on $\mathscr{M}$ having the standard extension property along the hypersurface $\{ f = 0 \}$. This gives a natural  $\mathcal{O}_{\mathscr{M}}$-linear embedding of this $\mathcal{O}_{\mathscr{M}}$-module in $Db_{\mathscr{M}}$.

\begin{cor}\label{no torsion}
The action of any  $P \in \mathcal{D}_{\mathscr{M}}$ on any element in the sub-$\mathcal{O}_{\mathscr{M}}$-module   $\mathcal{O}_{\mathscr{M}}\vert f\vert^{2\alpha}[ f^{-1}] \subset Db_{\mathscr{M}}$ coincides with the formal action defined above.\\
In particular, the sub-$\mathcal{D}_{\mathscr{M}}$-module generated by $Z_\alpha$ in $Db_{\mathscr{M}}$  has no torsion, where $Z_\alpha$ is the distribution on $\mathscr{M}$ associated to the locally bounded function $\vert f\vert^{2\alpha}$ on $\mathscr{M} $. 
\end{cor}

\parag{Proof} Let $V$ be a holomorphic vector field on $\mathscr{M} $. Then let $T_{\alpha,N}$ be the distribution defined by
$$ \langle T_{\alpha,N} , \xi \rangle := P_0(\lambda = \alpha, \int_\mathscr{M}  \vert f\vert^{2\alpha} f^{-N} \xi) $$
and let  $V^*$ be the adjoint of $V$;  we obtain:
\begin{align*}
& \langle V(T_{\alpha,N}) , \xi \rangle = \langle T_{\alpha,N} , V^*(\xi) \rangle = \lim_{\varepsilon \to 0} \int_{\vert f\vert \geq \varepsilon}  \vert f\vert^{2\alpha} f^{-N} V^*(\xi) \\ 
& \quad = \lim_{\varepsilon \to 0} \int_{\vert f\vert \geq \varepsilon} V(\vert f\vert^{2\alpha} f^{-N})\xi 
\end{align*}
because the boundary term on $\{ \vert f\vert = \varepsilon \}$ has limit $0$ when $\varepsilon$ goes to $0$ using the fact that 
 $$V(\vert f\vert^{2\alpha} f^{-N}) = (\alpha - N) V(f)  \vert f\vert^{2\alpha} f^{-N-1}, $$
 the same argument than in the proof of Theorem \ref{coincide 1} and the remark that $V(f)\xi$ is in $\mathscr{C}_c^\infty(\mathscr{M} )^{(n+1,n+1)}$.\\
 This gives, using again Theorem \ref{coincide 1}, that $V(T_{\alpha,N}) = V(f)T_{\alpha,N+1}$ in $Db_{\mathscr{M}}$.\\
 As holomorphic vector fields generate the $\mathcal{O}_\mathscr{M} $-algebra $\mathcal{D}_{\mathscr{M}}$, this is enough to complete the proof.$\hfill \blacksquare$\\
 
 \parag{Remark $R4$} The generalization of the previous corollary to the cases where we replace  $\vert f\vert^{2\alpha}$ by $\vert f\vert^{2\alpha}(Log \vert f\vert^2)^q$ for any positive integer $q$ is again an easy exercise.

\subsection{The case $\alpha  <  0$}

We shall explain now how to define, when $\alpha$ is a negative real number\footnote{We leave the case where $\alpha$ is a complex number with a negative real part as an exercise.} a  "Finite Part " of the integral 
$$ \int_{\mathscr{M}} \vert f\vert^{2\alpha}f^{-N}\xi$$
using the Asymptotic Expansion Theorem.
For  a differential form $\xi \in \mathscr{C}^\infty_c(\mathscr{M})^{(n+1, n+1)}$, we shall write as follows  the asymptotic expansion  when $s$ goes to $0$ of the fiber-integral 
 $$s \mapsto \int_{f=s} f\bar f \xi/df\wedge d\bar f  \simeq \sum_{r,j,m,m'} \tilde{T}_{m,m'}^{r,j}(\xi) \vert s\vert^{2r} s^m \bar s^{m'} (Log\vert s\vert)^j $$
 where $\xi \mapsto \tilde{T}(\xi) $ is a distribution with support in $\{f = 0 \}$ for each $r \in R$ a finite subset in $  [0, 1[ \cap \mathbb{Q}$, for $ j \in [0, n+1]$ and for $ m, m' \in \mathbb{N}$ (see the Corollary \ref{AS.EX.2} and also also  Remark 1 following Theorem \ref{AS.EX.1}).
 
 \begin{thm}\label{PV with pole}
 For $\alpha$ a negative real number such that $-\alpha \not\in R + \mathbb{N}$ the following limit exists and defines a distribution on $\mathscr{M}$ 
 \begin{align*}
 &  \lim_{\varepsilon \to 0}\Big( \int_{\{\vert f\vert \geq \varepsilon \}} \vert f\vert^{2\alpha} f^{-N}\xi  + 2i\pi \sum_{m = m'+N}^{ \alpha+r +m' < 0} \tilde{T}_{m,m'}^{r, j}(\xi) \frac{\varepsilon^{2(\alpha + r + m')}}{2(\alpha+r+m')}(-Log\varepsilon)^j\Big)
 \end{align*}
 which extends  to $\mathscr{M}$ the function $\vert f\vert^{2\alpha}f^{-N}$ on $\{ f \not= 0\}$.\\
 If $-\alpha$ is in $R + \mathbb{N}$, write $-\alpha = r + m', r \in R$ and $m' \in \mathbb{N}$ and add to  the sum  inside the limit in the left hand-side above the sum
 $$ 2i\pi\sum_{j = 0}^{n+1} \tilde{T}_{m'+N,m'}^{r, j}(\xi)(-Log\varepsilon)^{j+1}/(j+1) .$$
Then the result is analogous. 
 \end{thm}

 \parag{Proof} Using the asymptotic expansion for the fiber-integral $\int_{f=s} f\bar f \xi/df\wedge d\bar f $  recalled above, the proof that the limit in the left hand-side exists is analogous to the proof of Theorem \ref{PV.1}.$\hfill \blacksquare$ \\
 
  Of course, the constant term of the Laurent development at $\lambda = \alpha$  of the meromorphic extension of the  function $ \lambda \mapsto  \int_{\mathscr{M}} \vert f\vert^{2\lambda}f^{-N}\xi$ gives also such a distribution. 
  
 But in the case $\alpha < 0$ the relation between the distribution defined by  
 $$\xi \mapsto P_0\big(\lambda = \alpha, \int_\mathscr{M} \vert f \vert^{2\lambda}f^{-N}\xi\big)$$
  and the distribution defined in the previous theorem is not so clear in general, even if the difference between these two distributions has clearly its support in $\{ f = 0 \}$. \\
  Of course, when $\alpha$ satisfies the condition $\alpha + r > 0$ for each $r$ in $R \cup \{1\} \setminus \{0\}$, not only there is no term in the sum where $\alpha +r +m' \leq 0$ (using  Remark 1 following Theorem  \ref{AS.EX.1} for the case $0$) and the equality of the two distributions follows from the remark following Theorem \ref{coincide 1}.

\subsection{An easy generalization and an example}

Consider now $p \geq 2$ domains $\mathscr{M} _i$ in $\mathbb{C}^{n_i+1}$ for $i \in [1,p]$ and $p$ holomorphic functions $f_i : \mathscr{M} _i \to \mathbb{C}$. Then on $\mathscr{M}  := \prod_{i=1}^p \mathscr{M} _i$ it is easy  to prove, for $\alpha := (\alpha_1, \dots, \alpha_p)$ satisfying $\Re(\alpha_i) \geq 0, \forall i \in [1,p]$,  using Fubini' s Theorem, the existence of the distribution
$$ \langle T_{\alpha,N}, \xi \rangle := \lim_{\varepsilon \to 0} \langle T_\alpha^\varepsilon , \xi \rangle $$
where $\varepsilon := (\varepsilon_1, \dots, \varepsilon_p)$ is in $(\mathbb{R}^{+*})^p$, $N := (N_1, \dots, N_p) $ is in $ \mathbb{N}^p$, $\xi $ is in $\mathscr{C}^{\infty}_c(\mathscr{M} )^{(n+1, n+1)}$, $n+1 := \sum_{i=1}^p n_i+1$ and where
$$ \langle T_{\alpha,N}^\varepsilon , \xi \rangle := \int_{\prod_{i=1}^p \{\vert f_i\vert \geq \varepsilon_i \}} \quad \big( \prod_{i=1}^p \vert f_i\vert^{2\alpha_i}  f_i^{-N_i} \big)\xi .$$
It is also easy to make the meromorphic extension to $\lambda := (\lambda_1, \dots, \lambda_p) \in \mathbb{C}^p$ of the holomorphic distributions on $\mathscr{M} $ defined for $ \prod_{i=1}^p \{\Re(\lambda_i) \gg N_i \}$ and then to prove 
the  following generalization of Theorem \ref{coincide 1} to this "product" case:

\begin{thm}\label{generalized}
In the product situation described above the meromorphic extension of the distribution 
$$ \xi \mapsto \int_\mathscr{M}  \big(\prod_{i=1}^p \vert f_i\vert^{2\lambda_i}  f_i^{-N_i}\big) \xi $$
is holomorphic near the point $\alpha \in \mathbb{C}^p$ satisfying $\Re(\alpha_i) \geq 0, \forall i \in [1,p]$, and we have the equality
$$ P_0\Big(\lambda = \alpha, \int_\mathscr{M}  \big( \prod_{i=1}^p \vert f_i\vert^{2\lambda_i} f_i^{-N_i}\big) \xi\Big)  = \langle T_{\alpha,N}, \xi \rangle  $$
where the left hand-side denotes the value at $\lambda = \alpha$ of the meromorphic extension.
\end{thm}

\parag{Proof} The only point to precise in order to apply the $p = 1$ case and Fubini's Theorem successively to prove this generalization it the following remark:
\begin{itemize}
\item Let $\mathscr{M} _1\times \mathscr{M} _2 $ be the product of two complex manifolds $\mathscr{M} _1\times \mathscr{M} _2$ and let $\xi$ be a $\mathscr{C}_c^\infty $ test  differential form  on $\mathscr{M} _1\times \mathscr{M} _2$ and $T_2$ a distribution on $\mathscr{M} _2$. 
Then the test differential form defined on $\mathscr{M} _1$ by $\langle T_2, \xi \rangle $ is a  $\mathscr{C}_c^\infty $ test differential form on $\mathscr{M} _1$. So for any distribution $T_1$ on $\mathscr{M} _1$ the distribution $T_1\boxtimes T_2$ is well defined on $\mathscr{M} _1\times \mathscr{M} _2$ by the rule
 $$\langle T_1\boxtimes T_2, \xi \rangle := \langle T_1, \langle T_2, \xi \rangle \rangle.$$
\end{itemize}

\parag{Remark} Note that the definition of $T_{\alpha,N}$ implies that we have
$$ \langle T_{\alpha,N}, \xi \rangle = \lim_{\varepsilon \to 0} \int_{\vert \prod_{i=1}^p f_i \vert \geq  \varepsilon } \big(\prod_{i=1}^p \vert f_i\vert^{2\alpha_i}  f_i^{-N_i}\big)  \xi . $$
Then, using the argument in Corollary \ref{no torsion} on each $\mathscr{M} _i, i \in [1, p]$ we obtain that for any $P \in \mathcal{D}_{\mathscr{M}}$ the distribution $PT_{\alpha,N}$ has the standard extension property for the hypersurface $\{\prod_{i=1}^p f_i = 0 \}$. So the sub-$\mathcal{D}_{\mathscr{M}}$-module of $Db_{\mathscr{M}}$ generated by $T_{\alpha,N}$ has no torsion.\\
 Note also that this implies that the sub-$\mathcal{D}_{\mathscr{M}}$-module generated by $T_{\alpha,0}$  is contained in the $\mathcal{O}_\mathscr{M}$-module generated by the $T_{\alpha,N}$ when $N$ is in 
$\mathbb{N}^p$.

\parag{Remark $R5$} Again it is easy to generalize the previous theorem to the cases where we replace  $\vert f_i\vert^{2\alpha_i}$ by $\vert f_i\vert^{2\alpha_i}(Log \vert f_i\vert^2)^{q_i}$ for any non negative integers $q_i, i \in [1, p]$.

\parag{An example}

\parag{Notations} Let $\mathscr{M}  := \mathbb{C}^k$ with coordinates $z_1, \dots, z_k$ and let $\pi: \mathscr{M}  \to \mathscr{N}  \simeq \mathbb{C}^k$ be the quotient by the action of the permutation group $\mathfrak{S}_k$ on $\mathscr{M} $. Let $\sigma_1, \dots, \sigma_k$ be the elementary symmetric polynomials in $z_1, \dots, z_k$ which give a coordinate system on $\mathscr{N} $. Let also $\Delta := \prod_{1\leq i < j \leq k}  (z_i - z_j)^2 $ be the discriminant which is a polynomial in $\sigma_1, \dots, \sigma_k$.\\

\begin{lemma}\label{simplet}
Consider on $\mathscr{M} := \mathbb{C}^k$ with coordinates $z_1, \dots, z_k$ the holomorphic  functions $f := z_1$ and $\Delta := \prod_{1\leq i < j \leq k}  (z_i - z_j)^2 $. Then defining  new coordinates  $x_1 := z_1$ and  $x_h := z_h - z_1$ for $h \in [2, k]$ we obtain   a decomposition $\mathscr{M} = \mathscr{M}_1\times \mathscr{M}_2$ where $\mathscr{M}_1 := \mathbb{C}$ with coordinate $x_1$ and $\mathscr{M}_2 := \mathbb{C}^{k-1}$ with coordinates $x_h, h \in [2, k]$. Moreover  $f$ is in $\mathbb{C}[x_1]$ and $\Delta$ is in $\mathbb{C}[x_2, \dots, x_k]$.
\end{lemma}

\parag{Proof} The reader will see easily that  $\Delta = \prod_{h=2}^k x_h^2 \prod_{2\leq i < j \leq k} (x_i - x_j)^2$.$\hfill \blacksquare$\\

This easy lemma allows to apply Theorem \ref{generalized} to the situation described in the previous lemma.

\begin{cor}\label{no torsion M}
For any $\alpha \in \mathbb{C}$ such that $\Re(\alpha) \geq 0$ and any non negative integers $q$ and $N_1, N_2$ the distribution on $\mathscr{M}$  associated to the locally integrable function 
 $$X_{\alpha,N_1,N_2,q} := \big(\sum_{j=1}^k \vert z_j\vert^{2\alpha}(Log\vert z_j\vert^2)^q z_j^{-N_1} \big) \Delta^{-N_2}$$
 generates a sub-$\mathcal{D}_{\mathscr{M}}$-module of the $\mathcal{D}_{\mathscr{M}}$-module $Db_{\mathscr{M}}$ which has no torsion. Moreover, any distribution in this sub-$\mathcal{D}_{\mathscr{M}}$-module has the standard extension property along the divisor $\{\pi^{-1}(\sigma_k\Delta(\sigma) )= 0 \}$ in $\mathscr{M}$ and is $\mathscr{C}^\infty$(in fact real analytic)  outside this divisor.\\
 The same properties are true along the divisor $\{\sigma_k\Delta(\sigma) = 0\}$, for the sub-$\mathcal{D}_{\mathscr{N}}$-module of $Db_{\mathscr{N}}$ generated by the distribution $\pi_*(X_{\alpha,N_1,N_2,q})$.\\
\end{cor}

\parag{Proof}  For our assertion on $\mathscr{M}$ it is enough  to prove that for any holomorphic vector field $V$  on $\mathscr{M}$  the distribution $V(\vert z_1\vert^{2\alpha}z_1^{-N_1}\Delta^{-N_2})$ has the standard extension property and this is clear from the remark following Theorem \ref{generalized} and Remark $R5$. \\
The assertion on the corresponding  $\mathcal{D}_{\mathscr{N}}$-module  is consequence of the fact that if a distribution $T$ on $\mathscr{M}$ has the standard extension property along the divisor $\{\pi^{-1}(\sigma_k\Delta(\sigma)) = 0 \}$ then the distribution $\pi_*(T)$ has the standard extension property along the divisor $\{ (\sigma_k\Delta(\sigma) = 0 \}$ because the quotient map $\pi$ is a finite \'etale map outside $\{\Delta(\sigma) = 0\}$. $\hfill \blacksquare$

\parag{Remark} If $S$ and $T$ are two distributions on a complex manifold $\mathscr{N}$ such that $\mathcal{D}_{\mathscr{N}}S \subset Db_{\mathscr{N}}$ and $\mathcal{D}_{\mathscr{N}}T \subset Db_{\mathscr{N}}$ have no $\mathcal{O}_\mathscr{N}$-torsion, it is not clear that there is   no $\mathcal{O}_\mathscr{N}$-torsion in  $D_\mathscr{N}(S+T) \subset Db_{\mathscr{N}}$ . But if we know that $\mathcal{D}_{\mathscr{N}}S$ and $\mathcal{D}_{\mathscr{N}}T$ contain only distributions having the standard extension property for a given hypersurface $H$, then this also the case for any distribution
in $\mathcal{D}_{\mathscr{N}}S + \mathcal{D}_{\mathscr{N}}T \subset Db_{\mathscr{N}}$ which contains $\mathcal{D}_\mathscr{N}(S+T)$.

\section{Some conjugate D-modules}

The aim of this section is to give some examples of explicit computations of the conjugate module (in Kashiwara sense, see \cite{Ka1}) of some regular holonomic $\mathcal{D}$-modules. We consider here the case of the $\mathcal{D}$-modules associated to the multivalued functions $z(\sigma)^\lambda$ for $\lambda \in \mathbb{C}$ where $z(\sigma)$ is the root of the universal degree $k$ equation 
$$ z^k + \sum_{h=1}^k (-1)^h\sigma_hz^{k-h} $$
where $\sigma := (\sigma_1, \dots, \sigma_k)$ is in $N := \mathbb{C}^k$. The structure of these regular holonomic $\mathcal{D}_{\mathscr{N}}$-modules has been described in \cite{B4} and we describe here for each such  $\mathcal{D}$-module and  for each simple factor which appears in its decomposition a distribution $T$ on $\mathscr{N}$ which generates the sub-$\mathcal{D}_{\bar{\mathscr{N}}}$-module of $Db_{\mathscr{N}}$ which is the conjugate  of the module we consider.\\
The reader will notice that even if in each case the corresponding  distribution $T$ is rather easily constructed from the horizontal (multivalued) basis of the corresponding vector bundle with a simple pole connection associated to the $\mathcal{D}_{\mathscr{N}}$-module under consideration, the proof that this distribution generates the conjugate module uses in a crucial way the non trivial argument of "non torsion" which is proved in the previous section (see Corollary \ref{no torsion M}).\\

We begin by recalling the basic results on the conjugation functor of M. Kashiwara. The following theorem is proved in \cite{Ka1}:

\begin{thm}[Kashiwara conjugation functor]
Let $\mathscr{N}$ be a complex manifold, $\mathcal{D}_{\mathscr{N}}$ the sheaf of holomorphic partial differential operators on $\mathscr{N}$ and $\mathcal{D}_{\bar{\mathscr{N}}}$ the sheaf of  anti-holomorphic partial differential operator on $\mathscr{N}$. Note $Db_{\mathscr{N}}$ the sheaf of distributions on $\mathscr{N}$. It is a left-$\mathcal{D}_{\mathscr{N}}$-module but also a left-$\mathcal{D}_{\bar{\mathscr{N}}}$-module and these two actions commute.\\
For each regular holonomic $\mathcal{D}_{\mathscr{N}}$-module $\mathcal{N}$ the sub-$\mathcal{D}_{\bar{\mathscr{N}}}$-module 
$$c_\mathscr{N}(\mathcal{N}) := Hom_{\mathcal{D}_{\mathscr{N}}}(\mathcal{N}, Db_{\mathscr{N}})$$
 is regular holonomic (as a $\mathcal{D}_{\bar{\mathscr{N}}}$-module) and the contra-variant functor $c$ is an anti-equivalence of categories which satisfies $c_{\bar{\mathscr{N}}}\circ c_\mathscr{N} = Id$.$\hfill \blacksquare$
\end{thm}

Moreover, M. Kashiwara also obtains the following proposition which will be useful to describe the $\mathcal{D}_{\bar{\mathscr{N}}}$-modules $c(\mathcal{N}_\lambda)$:

\begin{prop}[ see \cite{Ka1} Prop. 5]\label{K2}
If $\mathcal{N}$ is a regular holonomic $\mathcal{D}_{\mathscr{N}}$-module on a complex manifold $\mathscr{N}$ and if $T$ is in $c(\mathcal{N})$, the following conditions are equivalent:
\begin{enumerate}[a)]
\item $T$ is an injective sheaf homomorphism of $\mathcal{N}$ to $Db_{\mathscr{N}}$.
\item  $T$ generates $c(\mathcal{N})$ as a $\mathcal{D}_{\bar{\mathscr{N}}}$-module.$\hfill \blacksquare$
\end{enumerate}
\end{prop}

\parag{Notations} We consider the quotient map  $ \pi : \mathscr{M} := \mathbb{C}^k \to  \mathscr{N}:= \mathbb{C}^k\big/\mathfrak{S}_k \simeq \mathbb{C}^k$ with respective coordinates $z_1, \dots, z_k$ and $\sigma_1,\dots, \sigma_k$, where $\sigma_h$ is the $h$-th  elementary symmetric function of $(z_1, \dots, z_k)$.\\
The vector fields on $\mathscr{M}$  associated to partial derivatives in $z_1, \dots, z_k$ are denoted by $\partial_{z_1}, \dots, \partial_{z_k}$ and the vector fields on $\mathscr{N}$  associated to partial derivatives in $\sigma_1, \dots, \sigma_k$ are denoted by $\partial_1, \dots, \partial_k$.\\
We note $\mathcal{I}$ the left ideal in $\mathcal{D}_{\mathscr{N}}$ generated by the following global sections:
\begin{align*}
&  A_{p,q} = \partial_p\partial_q - \partial_{p+1}\partial_{q-1} \quad {\rm where} \quad  (p,q) \in [1,k-1]\times [2,k] \\
& \mathcal{T}^m := \partial_1\partial_{m-1} + \partial_mE \quad {\rm where} \quad E := \sum_{h=1}^k  \sigma_h\partial_h \quad {\rm and} \quad m \in [2, k]
\end{align*}

Recall that a {\bf trace function} on $\mathscr{N}$ is a holomorphic function $F$ such that there exists a holomorphic function $f$ in one variable $z$ such that
$$ F(\sigma) = \sum_{j=1}^k f(z_j) $$
where $\sigma := ( \sigma_1, \dots, \sigma_k)$ are the elementary symmetric functions of $z_1, \dots, z_k$. \\
It is proved in \cite{B1} that the trace functions are annihilated by the left ideal $\mathcal{I}$ in $\mathcal{D}_{\mathscr{N}}$ and this characterizes the trace functions.\\

We note $U_{-1}$ the vector field on $\mathscr{N}$ given by $U_{-1} := k\partial_1 + \sum_{h=1}^{k-1} (k-h)\sigma_h\partial_{h+1}$ which is the image by the tangent map  $T_\pi$ to $\pi$  of the vector field  $V_{-1} := \sum_{j=1}^k \partial_{z_j}$ on $\mathscr{M}$.\\
We note $U_{0}$ the vector field on $\mathscr{N}$ given by  $U_0 := \sum_{h=1}^k  h\sigma_h\partial_h $ which is the image by $T_\pi$ of the vector field  $V_{0} := \sum_{j=1}^k z_j\partial_{z_j}$on $\mathscr{M}$.\\
We note $U_1$ the vector field on $\mathscr{N}$ given by $U_1 := \sum_{h=1}^k (\sigma_1\sigma_h -(h+1)\sigma_{h+1})\partial_h$ which is the image by $T_\pi$ of the vector field  $V_{1} := \sum_{j=1}^k z_j^2\partial_{z_j}$ on $\mathscr{M}$.\\
The left  ideal $\mathcal{J}_\lambda$ in $\mathcal{D}_{\mathscr{N}}$ is, by definition, the sum  $\mathcal{I} + \mathcal{D}_{\mathscr{N}}(U_0 - \lambda)$ and we define  $\mathcal{N}_\lambda := \mathcal{D}\big/\mathcal{J}_\lambda$.\\

The following results which give the structure of the $\mathcal{D}_{\mathscr{N}}$-module   $\mathcal{N}_\lambda$ for each $\lambda \in \mathbb{C}$  are proved in \cite{B4}:\\
\begin{enumerate}
\item For  each complex number $\lambda$, the  $\mathcal{D}_{\mathscr{N}}$-module   $\mathcal{N}_\lambda$  is   holonomic and regular.
\item The right multiplication by $U_{-1}$ induces a  $\mathcal{D}_{\mathscr{N}}$-linear map  
$$ \square U_{-1}  : \mathcal{N}_\lambda \to \mathcal{N}_{\lambda+1}$$
which is  an isomorphism for each $\lambda \not= -1, 0$. Moreover the right multiplication by $U_1$ induces an isomorphism $ \square U_1 : \mathcal{N}_{\lambda +1} \to \mathcal{N}_{\lambda }$ for any $\lambda \not= 0, -1$ and we have $\square U_1\circ \square U_{-1} = \lambda(\lambda+1)$ on $\mathcal{N}_{\lambda}$.
\item For $\lambda \not\in \mathbb{Z}$ the $\mathcal{D}_{\mathscr{N}}$-module $\mathcal{N}_\lambda$ is simple.
\item The kernel  $\mathcal{N}_{-1}^\square$ of the $\mathcal{D}_{\mathscr{N}}$-linear map $\varphi_{-1} : \mathcal{N}_{-1} \to \mathcal{O}_N(\star \sigma_k)$ defined by $\varphi_{-1}(1) = \sigma_{k-1}/\sigma_k$ is simple.
\item The sub-module $\mathcal{N}_0^\square$ generated by $U_1$ in $\mathcal{N}_0$ is simple and the quotient $\mathcal{N}_0\big/\mathcal{N}_0^\square$ is isomorphic to $\mathcal{O}_N(\star \sigma_k)$ 
\item The torsion sub-module $\mathfrak{T}$ in $\mathcal{N}_1$ is isomorphic to $H^1_{[\sigma_k = 0]}(\mathcal{O}_N)$; it  is generated by the class of $\partial_kU_{-1}$ in $\mathcal{N}_1$. 
\item The sub-module $Im(\square U_{-1})$in  $\mathcal{N}_1$ (which is generated by $U_{-1}$)  contains $\mathfrak{T}$ and the quotient $Im(\square U_{-1})\big/\mathfrak{T}$ is isomorphic to $\mathcal{O}_\mathscr{N}$ via the map $\varphi_1 :\mathcal{N}_1 \to \mathcal{O}_N$ defined by $\varphi_1(1) = \sigma_1$   (and then $[U_{-1}] \mapsto k$).
\item The quotient $\mathcal{N}_1^\square := \mathcal{N}_1\big/Im(\square U_{-1})$ is simple and isomorphic to the quotient
$$\mathcal{D}_{\mathscr{N}}\big/ \mathcal{I} + \mathcal{D}_{\mathscr{N}}(U_0-1) + \mathcal{D}_{\mathscr{N}}U_{-1}.$$
\item The right multiplication by $U_1$ which sends $\mathcal{N}_1$ in $\mathcal{N}_0^\square$ vanishes\footnote{It is proved in \cite{B4} Formula $(19)$ page 20,  that we have $U_{-1}U_1 = (U_0+1)U_0 \quad {\rm modulo} \quad \mathcal{I}$.} on $Im(\square U_{-1})$ and  induces an isomorphism of  $\mathcal{N}_1^\square$ onto $\mathcal{N}_0^\square$.
\item The right multiplication by $U_{-1}$ which sends $\mathcal{N}_{-1}$ to  $\mathcal{N}_0 $   induces an isomorphism of  $\mathcal{N}_{-1}^\square$ onto $\mathcal{N}_0^\square$.
\end{enumerate}

\parag{Important remark} The point 2  recalled above shows that to study the $\mathcal{D}_{\mathscr{N}}$-modules $\mathcal{N}_\lambda$ it is enough to consider the following cases:

\begin{itemize}
\item The cases $\Re(\lambda) \in [0, 1[$ and $\lambda \not= 0$. We call it the case $G$. 
\item The cases $\lambda = -1, 0, 1$. We call them the case $\lambda = -1, 0, 1$ respectively.
\end{itemize}
Then for any $\lambda \in \mathbb{C}$ we reach one of these previous cases using an isomorphism given either by $U_1^N$ or by $U_{-1}^N$ for a suitable $N \in \mathbb{N}$.\\

Now define the following distributions on $\mathscr{M}$:
\begin{enumerate}
\item For $\Re(\lambda) \in [0, 1[, \lambda \not= 0$ define $X_\lambda := \sum_{j=1}^k \vert z_j\vert^{2\lambda} $. 
\item  $X_1 := \sum_{j=1}^k  \vert z_j - \sigma_1/k\vert^2 =  \sum_{j=1}^k  \vert z_j\vert^2 -\vert\sigma_1\vert^2/k$.
\item $X_0 := \sum_{j=1}^k (\bar z_j-\bar \sigma_1/k)Log \vert z_j\vert^2  = \sum_{j=1}^k  \vert z_j\vert^2z_j^{-1}Log \vert z_j\vert^2  - (\bar \sigma_1/k) Log\vert \sigma_k\vert^2$.
\item $X_{-1} = \sum_{j=1}^k (\bar z_j- \bar\sigma_1/k) z_j^{-1}  =   \sum_{j=1}^k  \vert z_j \vert^2z_j^{-2} - (\bar\sigma_1/k)\sigma_{k-1}/\sigma_k$.

\end{enumerate}
For each case we define the distribution $\mathcal{X}_\lambda := \pi_*(X_\lambda)$ on 
$\mathscr{N}$.\\

Define also  the following distributions on $\mathscr{N}$:

\begin{enumerate}
\item For the case $G$ :\quad  $\mathcal{Y}_\lambda := 0$.
\item  For the case $\lambda = 1$:  \quad $\mathcal{Y}_1 := \sigma_1/\bar\sigma_k$ 
\item  For the case $\lambda = 0 $: \quad  $\mathcal{Y}_0 := 1/\bar \sigma_k$.
\item  For the case $\lambda = -1$:  \quad $\mathcal{Y}_{-1} := \sigma_{k-1}/\sigma_k$.
\end{enumerate}


Then we have the following results:

\begin{thm}\label{conj.1}[conjugate modules]\\
In all cases $G, -1, 0, 1$ the distribution $\mathcal{X}_\lambda$ defines an element in $c(\mathcal{N}_\lambda)$ via the $\mathcal{D}_{\mathscr{N}}$-linear map sending $1$ to $\mathcal{X}_\lambda$.\\
 In case $G$  the distribution $\mathcal{X}_\lambda$ defines a   generator  of the (simple) $\mathcal{D}_{\bar{\mathscr{N}}}$-module $c(\mathcal{N}_\lambda)$.\\
In cases $\lambda =1$  the distribution $\mathcal{X}_1$ defines an element  in $c(\mathcal{N}_1^\square) \subset c(\mathcal{N}_1)$  and  gives a generator of this (simple) $\mathcal{D}_{\bar{ \mathscr{N}}}$sub-module.\\
In cases $\lambda = -1, 0$ the simple module $c(\mathcal{N}_\lambda^\square)$ is a quotient of $c(\mathcal{N}_\lambda)$ and the image of $\mathcal{X}_\lambda $ in  $c(\mathcal{N}_\lambda^\square)$ gives a generator of this (simple) $\mathcal{D}_{\bar{ \mathscr{N}}}$-module.\\
In all cases $\lambda = -1, 0, 1$ the map sending $1$ to $\mathcal{Y}_\lambda$ is in $c(\mathcal{N}_\lambda)$ and the $\mathcal{D}_{\bar{\mathscr{N}}}$-module  $c(\mathcal{N}_\lambda)$  is generated by  $\mathcal{X}_\lambda$ and $\mathcal{Y}_\lambda$. 
\end{thm}

\parag{Proof} Consider the first point. We have to prove that the left ideal  $\mathcal{J}_\lambda$ annihilates the distribution $\mathcal{X}_\lambda$ for each  case.\\
Corollary \ref{no torsion M} and Remark R5  give that for $\Re(\lambda) \geq 0$ any distribution in  the $\mathcal{D}_\mathscr{M}$-module generated by the distribution 
$$ \vert z_j\vert^{2\lambda}(Log \vert z_j\vert^2)^q z_j^{-N}$$
has the standard extension property. So the same property holds for any distribution in the $\mathcal{D}_\mathscr{N}$-module generated by the distribution 
$$ \sum_{j=1}^k  \vert z_j\vert^{2\lambda}(Log \vert z_j\vert^2)^q z_j^{-N}$$
along the hypersurface $\{\sigma_k\Delta(\sigma) = 0 \}$  thanks to  Corollary \ref{no torsion M}  .\\
Remark also that the generator of the left ideal $\mathcal{I}$ and the vector field $U_0$ are in the left ideal of $\mathcal{D}_\mathscr{N}$ generated by $\partial_h, h \in [1, k-1]$ and $\sigma_k\partial_k$ and that these vector fields annihilate the distribution $1/\bar \sigma_k$ (see Lemma \ref{facile} below).\\
We have also  $\sigma_k\partial_k(Log \vert\sigma_k\vert^2) = 1$ in $Db_{\mathscr{N}}$. \\
Now the verification that in each case the distribution $\mathcal{X}_\lambda$ is annihilated by the generator of the left ideal $\mathcal{J}_\lambda$ is "formal"  (see Corollary \ref{no torsion M} ) because $\mathcal{X}_\lambda$ is locally on the open set  $\{ \sigma_k\Delta(\sigma) \not= 0\}$ a linear combination with anti-holomorphic coefficients of the holomorphic trace functions $z_j(\sigma)^\lambda, j \in [1, k]$.\\
Moreover it is clear that $z_j(\sigma)^\lambda$ is homogeneous of degree $\lambda$ on $\mathscr{M}$ and the vector field   $U_0$ is equal to $ \pi_*(V_0)$ where $V_0$  is the Euler vector field on $\mathscr{M}$. \\

  In the case $G$ we know that $c(\mathcal{N}_\lambda)$ is a simple  $\mathcal{D}_{\bar{\mathscr{N}}}$-module. So the second point of the theorem is proved as $\mathcal{X}_\lambda$ is clearly not $0$ in $Db_\mathscr{N}$.\\
 
For $\lambda = 1$, as we know that ideal  $\mathcal{J}_1 = \mathcal{I} + \mathcal{D}_{\mathscr{N}}(U_0-1)$ already annihilates $\mathcal{X}_1$ it is enough, thanks to the point 8 recalled above, to check that $U_{-1}(\mathcal{X}_1) = 0$ in $Db_{\mathscr{N}}$ to conclude because we know that $\mathcal{N}_1^\square$ is simple and $\mathcal{X}_1$ is not $0$ in $Db_{\mathscr{N}}$. Again the fact that it is enough to check this on the open set $\{\Delta(\sigma) \not= 0\}$ makes the computation easy :
$$ V_{-1}(\vert z_j\vert^2) = \bar z_j \quad {\rm and} \quad V_{-1}(\vert \sigma_1\vert^2) = k\bar \sigma_1 \quad {\rm so} \quad U_{-1}(\mathcal{X}_1) = 0.$$

\begin{lemma}\label{facile}
Let $\mathcal{M}$ be the $\mathcal{D}_{\mathbb{C}}$-module quotient of $\mathcal{D}_{\mathbb{C}}$ by the left ideal generated by $z\partial_z$. Then $c(\mathcal{M})$ is isomorphic to $\mathcal{O}_{\bar{\mathbb{C}}}(*\bar z)$.
\end{lemma}

\parag{Proof} First recall that the distribution on $\mathbb{C}$ associated to the locally integrable function $1/\bar z$ satisfies $\partial_z(1/\bar z) = i\pi \delta_0$, where $\delta_0$ is the Dirac mass at the origin.
So the $\mathcal{D}_{\overline{\mathbb{C}}}$-sub-module of $Db_{\mathbb{C}}$ generated by $1/\bar z$ is contained in  $c(\mathcal{M})$. To prove that it is equal to $c(\mathcal{M})$, thanks to Proposition \ref{K2}, it is enough to prove that the kernel of the $\mathcal{D}_{\mathbb{C}}$-linear map $\mathcal{M} \to Db_\mathbb{C} $ defined by $1 \mapsto 1/\bar z$ is injective. But modulo the left ideal generated by $z\partial_z$ any 
$P \in \mathcal{D}_{\mathbb{C}}$ may be written as $a_0 + \sum_{p=1}^m a_p\partial_z^p$ where $a_0, \dots, a_p$ are complex numbers. Then $a_0$ has to vanish and  the $\mathbb{C}$-linear independence of the distributions $\partial_z^{p-1}\delta_0, p \geq 1$ allows to conclude.$\hfill \blacksquare$

\begin{cor}\label{utile}
The conjugate of the sub-$\mathcal{D}_{\mathscr{N}}$-module $Im(\square U_{-1}) := \mathcal{D}_{\mathscr{N}}U_{-1} \subset \mathcal{N}_1$ is equal to   $\mathcal{D}_{\bar{\mathscr{N}}}\sigma_1/\bar\sigma_k \subset Db_{\mathscr{N}}$ which is isomorphic to the $\mathcal{D}_{\bar{\mathscr{N}}}$-module  $\mathcal{O}_{\bar{\mathscr{N}}}(*\bar \sigma_k)$ . It is generated by the $\mathcal{D}_{\mathscr{N}}$-linear map sending the class of $U_{-1}$ in $Im(\square U_{-1}) \subset \mathcal{N}_1$ to the distribution $\mathcal{Y}_1 = \sigma_1/\bar \sigma_k$. So the distribution $\mathcal{X}_1$ and $ \mathcal{Y}_1$ in $Db_{\mathscr{N}}$  generate $c(\mathcal{N}_1)$ as a  $\mathcal{D}_{\bar{\mathscr{N}}}$-sub-module of $Db_{\mathscr{N}}$.
\end{cor}

\parag{Proof} It is proved in Proposition 4.1.7 and Lemma 4.1.8 of \cite{B4} that the annihilator of the class of $U_{-1}$ in $\mathcal{N}_1$ is generated by $\partial_h, h \in [1, k-1]$ and $\sigma_k\partial_k$. Note that this is also the annihilator in $\mathcal{D}_{\mathscr{N}}$ of the distribution $1/\bar \sigma_k$ in $Db_{\mathscr{N}}$ thanks to the previous  lemma. This shows that the $\mathcal{D}_{\mathscr{N}}$-linear map defined by  $[U_{-1}] \mapsto 1/\bar\sigma_k \in Db_{\mathscr{N}}$ is injective, proving the isomorphism of $\mathcal{D}_{\bar{\mathscr{N}}}$-modules  $c(Im(\square U_{-1})) \simeq \mathcal{D}_{\bar{\mathscr{N}}}1/\bar\sigma_k \simeq \mathcal{O}_{\bar{\mathscr{N}}}(*\bar \sigma_k) $ thanks to Proposition \ref{K2}.\\
Then the exact sequence of $\mathcal{D}_{\mathscr{N}}$-modules
$$ 0 \to Im(\square U_{-1}) \to \mathcal{N}_1 \to \mathcal{N}_1^\square \to 0 $$
gives the exact sequence of $\mathcal{D}_{\bar{\mathscr{N}}}$-modules
$$ 0 \to c(\mathcal{N}_1^\square) \to c(\mathcal{N}_1) \overset{\alpha}{\to} 
c(Im(\square U_{-1}) \to 0 $$
where the map $\alpha$ is the quotient map. Now, let us prove that  the distribution 
 $\mathcal{Y}_1 = \sigma_1/\bar \sigma_k$ is  a generator of  $c(Im(\square U_{-1}))$ via the $\mathcal{D}_{\mathscr{N}}$-linear map   $U_{-1} \mapsto \mathcal{Y}_1$.
 As $\sigma_1$ commutes with the action of $\mathcal{D}_{\bar{\mathscr{N}}}$, the annihilator of $\sigma_1/\bar \sigma_k$ in $\mathcal{D}_{\bar{\mathscr{N}}}$ is the same than the annihilator of $1/\bar \sigma_k$ proving our claim. \\
But  we know that $\mathcal{X}_1$ is a generator of $c(\mathcal{N}_1^\square)$.  So, it is enough to prove that $\mathcal{Y}_1$ is in $c(\mathcal{N}_1)$ to see that with $\mathcal{X}_1$ they generate this $\mathcal{D}_{\bar{\mathscr{N}}}$-module. This is proved in our next lemma and completes the proof of this corollary.$\hfill \blacksquare$\\

\begin{lemma}\label{ajout}
For $\lambda = 1, 0, -1$ the distribution  $\mathcal{Y}_\lambda$ defines an element in $c(\mathcal{N}_\lambda)$ via the $\mathcal{D}_{\mathscr{N}}$-linear map $1 \mapsto \mathcal{Y}_\lambda$.\\
\end{lemma}

\parag{Proof} The fact that $\mathcal{Y}_1$ in in $c(\mathcal{N}_1)$ is easy because $\sigma_1$ is a trace function which satisfies $U_0(\sigma_1) = \sigma_1$, because $\partial_h, h \in [1, k-1]$ and 
$ \sigma_k\partial_k$ annihilate the distribution $1/\bar\sigma_k$ and it is easy to check that the generators of $\mathcal{I}$ and $U_0$ are in this left ideal of $\mathcal{D}_{\mathscr{N}}$ generated by these vector fields.\\
Now we have also, using the same argument, that $\mathcal{Y}_0$ is in $c(\mathcal{N}_0)$.\\
To see that $\mathcal{Y}_{-1}$ is annihilated by $\mathcal{J}_{-1}$, remark first that $\mathcal{D}_{\mathscr{N}}\mathcal{Y}_{-1} \subset Db_{\mathscr{N}}$ has no torsion (see Corollary \ref{no torsion M} for $\alpha = 0, N = -1, m = 0$).
 Then note that
$$\sigma_{k-1}/\sigma_k = \pi_*(\sum_{j=1}^k 1/z_j)$$
 is a trace function on the open set  $\{\sigma_k\Delta(\sigma) \not= 0 \}$ and that $U_0 + 1$ kills also this holomorphic function on this open set. This is enough to complete the proof of the lemma using again the absence of torsion.$\hfill \blacksquare$\\
 
 \begin{lemma}\label{autre ajout}
The $\mathcal{D}_{\mathscr{N}}$-linear map defined by  $[1] \mapsto \mathcal{X}_0 \in Db_{\mathscr{N}}$  is a generator of $c(\mathcal{N}_0^\square)$ and the $\mathcal{D}_{\mathscr{N}}$-linear map defined by  $[1] \mapsto \mathcal{X}_{-1} \in Db_{\mathscr{N}}$  is a generator of $c(\mathcal{N}_{-1}^\square)$.
\end{lemma}

\parag{Proof} As we already know that  the distribution $\mathcal{X}_1 $ is a generator of $c(\mathcal{N}_1^\square)$ the point 9 above implies that a distribution $Z$  in $c(\mathcal{N}_0)$ which satisfies $U_1(Z) = \mathcal{X}_1$ must induce a generator in $c(\mathcal{N}_0^\square)$. But $\mathcal{X}_0$ is such a distribution.\\
Thanks to point 10 above  and the previous result, the distribution $U_{-1}(\mathcal{X}_0)$  gives a generator in $c(\mathcal{N}_{-1}^\square)$ and we have, using as above the absence of torsion the equality $U_{-1}(\mathcal{X}_0) = \mathcal{X}_{-1}$ in $Db_{\mathscr{N}}$.$\hfill \blacksquare$\\

\parag{End of proof of Theorem \ref{conj.1}} For $\lambda = 0$ we have the exact sequence of $\mathcal{D}_{\mathscr{N}}$-modules (see point 5 above)
$$ 0 \to \mathcal{N}_{0}^\square \to \mathcal{N}_0 \to \mathcal{O}_N(*\sigma_k) \to 0 $$
where the simple sub-module $\mathcal{N}_0^\square$ is generated by $U_1$. So we have an exact sequence of $\mathcal{D}_{\bar{\mathscr{N}}}$-modules
$$ 0 \to c( \mathcal{O}_N(*\sigma_k)) \to c(\mathcal{N}_0) \to c(\mathcal{N}_0^\square) \to 0 $$
and if we find a distribution $Z$ in $c(\mathcal{N}_0) $ such that its image in $c(\mathcal{N}_0^\square)$ is equal to $\mathcal{X}_{-1}$, and a distribution $T$ which generates $c( \mathcal{O}_N(*\sigma_k)) $,  then $Z$ and $ T$ will  generate  $c(\mathcal{N}_0)$. \\
But $U_1(\mathcal{Y}_{0}) = 0$ in $Db_{\mathscr{N}}$ because the coefficient of $\partial_k$ in $U_1$ is $\sigma_1\sigma_k$ and $\sigma_k\partial_k$ kills the distribution $1/\bar\sigma_k$. So $\mathcal{X}_0$ and $ \mathcal{Y}_0$ generates $c(\mathcal{N}_0)$.\\
The situation is analogous for $\lambda = -1$ because $\mathcal{X}_{-1}$ is in $c(\mathcal{N}_{-1})$ and its image generates $c(\mathcal{N}_{-1}^\square)$ thanks to the previous lemma. \\
Now the distribution $\mathcal{Y}_{-1}$ is in $c(\mathcal{N}_{-1})$ and generates a sub-$\mathcal{D}_{\bar N}$-module isomorphic to $\mathcal{O}_{\bar N}(*\bar \sigma_k)$ because its annihilator is generated by $\bar\partial_h, h \in [1, k-1]$ and $\bar\sigma_k\bar\partial_k$. Its image in $c(\mathcal{N}_{-1}^\square)$ is either a generator or $0$. But it has to be non $0$ because there is no non zero morphism between $\mathcal{O}_{\bar N}(*\bar \sigma_k)$ and $c(\mathcal{N}_{-1}^\square)$. So $\mathcal{X}_{-1}$ and $ \mathcal{Y}_{-1}$ generate $c(\mathcal{N}_{-1})$.$\hfill \blacksquare$\\

\parag{Complement}

We  give in the following lemma, for each case $\lambda = -1, 0, 1$, a distribution which generates the $\mathcal{D}_{\bar{\mathscr{N}}}$-module $c(\mathcal{N}_\lambda)$.

\begin{lemma}
Consider the distribution $Z_{-1} := \sum_{j=1}^k \bar z_j/z_j$ on $\mathscr{M}$ and define \\ $\mathcal{Z}_{-1} := \pi_*(Z_{-1})$. Then $\mathcal{Z}_{-1}$ is a generator of $c(\mathcal{N}_{-1})$.
\end{lemma}

\parag{Proof} First write $Z_{-1} = \sum_{j=1}^k \vert z_j\vert^2/z_j^2 $ in order to show that  any distribution in the  $\mathcal{D}_{\mathcal{M}}$-module generated by $Z_{-1}$ has the standard extension property,
by applying Corollary \ref{no torsion M}. Then the same is true for any distribution in the 
$\mathcal{D}_{\mathcal{N}}$-module generated by $\mathcal{Z}_{-1}$.
Then it is easy to see that $\mathcal{Z}_{-1}$ is in $c(\mathcal{N}_{-1})$.\\
Then we have  $\bar U_{-1}(\mathcal{Z}_{-1}) = \pi_*(\sum_{j=1}^k 1/z_j )= \mathcal{Y}_{-1}$ and 
$$ \mathcal{X}_{-1} = \mathcal{Z}_{-1} - (\bar\sigma_1/k) \sigma_{k-1}/\sigma_k = (1 -( \bar \sigma_1/k)\bar U_{-1})( \mathcal{Z}_{-1}) \in  \mathcal{D}_{\bar{\mathscr{N}}} \mathcal{Z}_{-1}.$$
So $\mathcal{D}_{\bar{\mathscr{N}}} \mathcal{Z}_{-1}$ is equal to $c(\mathcal{N}_{-1})$ thanks to Theorem \ref{conj.1}.$\hfill \blacksquare$\\

\begin{lemma}
Let $Z_1 := \sum_{j=1}^k \vert z_j \vert^2 - \sigma_1/\sigma_k$ in $Db_{\mathscr{M}}$ and $\mathcal{Z}_1 := \pi_*(Z_1) \in Db_{\mathscr{N}}$. Then $\mathcal{Z}_{1}$ is a generator of $c(\mathcal{N}_{1})$.
\end{lemma}

\parag{Proof}  Any distribution in $\mathcal{D}_{\mathscr{N}} \mathcal{Z}_1$ has the standard extension property because this is true for the $\mathcal{D}_{\mathscr{N}}$-module generated by  $\pi_*(\sum_{j=1}^k \vert z_j \vert^2 )$ and and by $1/\sigma_k$.
So $\mathcal{Z}_1 $ is in $c(\mathcal{N}_1)$. Moreover we have
$$ (\bar U_0 +k)(\mathcal{Z}_1 ) = (k+1)\pi_*(\sum_{j=1}^k \vert z_j \vert^2 )= (k+1) \mathcal{X}_1 + \frac{k+1}{k}\bar\sigma_1\bar\sigma_k \mathcal{Y}_1 $$
and
$$ (\bar U_0 -1)(\mathcal{Z}_1 ) = (k-1)\mathcal{Y}_1.$$
So $\mathcal{D}_{\bar{\mathscr{N}}} \mathcal{Z}_{1}$ is equal to $c(\mathcal{N}_{1})$ thanks to Theorem \ref{conj.1}.$\hfill \blacksquare$\\

\begin{lemma}
 $\mathcal{X}_{0} + \mathcal{Y}_0 $ is a generator of $c(\mathcal{N}_{0})$.
\end{lemma}

\parag{Proof} We already know that $\mathcal{X}_0$  and $\mathcal{Y}_0$ are in $c(\mathcal{N}_{0})$.  Moreover we have
$$ (\bar U_0 -1)(\mathcal{X}_0 +  \mathcal{Y}_0) = -(k+1)/\bar \sigma_k =  -(k+1)\mathcal{Y}_0 $$
and
$$ (\bar U_0 + k )(\mathcal{X}_0 + 1/\bar \sigma_k) = \mathcal{X}_0  .$$
This is enough to conclude, thanks to Theorem \ref{conj.1}.$\hfill \blacksquare$\\


\providecommand{\bysame}{\leavevmode\hbox to3em{\hrulefill}\thinspace}
\providecommand{\MR}{\relax\ifhmode\unskip\space\fi MR }
\providecommand{\MRhref}[2]{%
  \href{http://www.ams.org/mathscinet-getitem?mr=#1}{#2}
}
\providecommand{\href}[2]{#2}
\begin{thebibliography}{}

\end{thebibliography}


\newpage



\begin{thebibliography}{99}

\parag{Reference}

\bibitem{B1} D. Barlet  {\em On symmetric partial differential operators}, math.arXiv:1911.09347,  to appear in Math. Zeitschrift.
\bibitem{B2} D. Barlet  {\em D\'eveloppement asymptotique des fonctions obtenues par int\'egration sur les fibres},
Invent. Math. \textbf{68}, 1, 129-174 (1982).
\bibitem{B3} D. Barlet  {\em  Fonctions de type trace}, Ann. Inst. Fourier \textbf{33}, 2, 43-76 (1983).
\bibitem{B4} D. Barlet {\em On partial differential operators which annihilate the roots of the universal equation of degree k}, math.arXiv:2101.01895 

\bibitem{B-K} Barlet, D. et Kashiwara, M. {\it Le r\'eseau $L^2$ d'un syst\`{e}me holonome r\'egulier} Inv. math. vol. 83 (1986) pp. 35-62.

\bibitem{BMa} D. Barlet and  H.-M. Maire  {\em Asymptotique des int\'egrales-fibres}, Ann. Inst. Fourier  \textbf{43}, 5, 1267-1299 (1993).


\bibitem{B-M}
D. Barlet and J. Magn\'usson
{\em Cycles analytiques complexes I: th\' eor\`emes de pr\'eparation des cycles,}
Cours Sp\'ecialis\'es  \textbf{22}, Soci\'et\'e Math\'ematique de France  (2014).

\bibitem{Be}
  Bernstein, J. N.  {\it Prolongement analytique des fonctions g\'en\'eralis\'ees avec param\`{e}tres} (en russe) Funkts Analyz 6.4 (1972).

\bibitem{Bj1} Bj\"ork, J. E. {\it Ring of differential operators} North Holland (1979)

\bibitem{Bj2}
 Bj\"ork, J. E.  {\it Analytical $D$-Modules and Applications,}
Mathematics and Its Applications, Kluwer Academic Publishers \textbf{247} (1993).


\bibitem{H-L} Herrera,M. and Lieberman, {\it Residues and principal values on complex spaces}, Math. Annalen 194, 259-294 (1971).



\bibitem{K} 
Kashiwara, M. {\it b-function and holonomic systems, rationality of roots of b-functions}, Inv. Math. 38 (1976) pp. 33-53.

\bibitem{Ka1} M. Kashiwara,
{\em Regular Holonomic $D$-modules and Distributions on a Complex Manifolds,}
\/ Advanced Studies in Pure Mathematics {\bf 8}, 1986, Complex Analytic Singularities pp. 199-206.














\end{thebibliography}
\end{document}